\documentclass{amsart}
\usepackage{amssymb}
\usepackage{amsthm}
\usepackage{amsmath}
\usepackage{url,float}
\usepackage{graphicx}
\usepackage{enumitem}
\usepackage[colorlinks]{hyperref}
\usepackage{tikz-cd}
\usepackage{geometry}
\usepackage{bbold}
\usepackage{xcolor}
\usepackage{todonotes}
\usepackage{listings}
\usepackage{booktabs}
\usepackage{makecell}
\usepackage[most]{tcolorbox}

\newcommand{\F}{\mathbb{F}}
\newtheorem{thm}{Theorem}[section]
\newtheorem{cor}[thm]{Corollary}
\newtheorem{lem}[thm]{Lemma}

\newtheorem{defn}[thm]{Definition}

\newtheorem{example}[thm]{Example}
\newcommand{\Z}{\mathbb{Z}}

\title[On Lattice Tilings by Symmetric Limited Magnitude Error Balls]{On Lattice Tilings by Symmetric Limited Magnitude \\Error Balls: Q-Polynomials and Ideal Characterization}
\author[A.\ Xiao]{Ang Xiao}
\email{xiaoang\_nudt@163.com}
\address{College of Science, National University of Defense Technology, 410073 Changsha, China}

\begin{document}
	\begin{abstract}
		In this paper, we study the lattice tilings of $\Z^n$ by the symmetric limited magnitude error balls $\mathcal{B}(n,t,k,k)$, which play a crucial role in flash memory applications. 
		Using generating functions and complete Bell polynomials, we derive an explicit formula for the $r$-th Q-polynomial associated with $\mathcal{B}(n,t,k,k)$. 
		Then, an ideal-theoretic characterization of the lattice tiling condition is developed: for an odd prime $p$ dividing $|\mathcal{B}(n,t,k,k)|$ exactly once, a lattice tiling would force a certain polynomial to lie outside a certain ideal of $\F_p[Y_1,\dots,Y_{\frac{p-1}{2}}]$. 
		Applying this criterion yields congruence conditions on the dimension $n$ that rule out lattice tilings; for $2\le t\le10$ and $1\le k\le10$, a systematic Maple computation over all primes $3\le p\le29$ excludes, on average, $35.65\%$ of all dimensions, with a maximum of $76.69\%$. 
		Finally, we obtain two infinite families of parameters $(n,t,k)$ for which no lattice tiling exists: no lattice tiling exists if $k\equiv4\pmod9$ and the base-$3$ addition of $t$ and $n-1-t$ produces exactly one carry; and no lattice tiling of $\Z^n$ by $\mathcal{B}(n,n-1,k,k)$ exists if $k\equiv8\pmod9$ and $n\equiv2,4\pmod6$. Their sizes are also estimated.
	\end{abstract}
	\subjclass[2020]{Primary 52C22; Secondary 05E05, 94B25, 11T06}
	\keywords{lattice tiling, limited magnitude error ball, perfect code}
	\maketitle
	\section{Introduction}
	For integers $n\ge t\ge 1$ and $k_1,k_2\ge 0$, the \emph{limited magnitude error ball} is defined as 
	\begin{align*}
		\mathcal{B}(n,t,k_1,k_2):=\{\mathbf{a}=(a_1,\dots,a_n)\in\Z^n:a_i\in[-k_2,k_1],\ \text{wt}(\mathbf{a})\le t\},
	\end{align*}
	where $\text{wt}(\mathbf{a})$ denotes the Hamming weight.  A subset $C\subset\Z^n$ forms a \emph{tiling} by $\mathcal{B}(n,t,k_1,k_2)$ if the translates $\{\mathcal{B}(n,t,k_1,k_2)+c:c\in C\}$ partition $\Z^n$; when $C$ is a lattice, we speak of a \emph{lattice tiling}.  Such tilings correspond to linear perfect codes under the limited magnitude error model, which is particularly relevant to flash memory systems \cite{cassuto_codes_2010}.
	
	\subsection*{Previous work}
	The case $t=1$ has been extensively studied.  The \emph{cross} $\mathcal{B}(n,1,k,k)$ and \emph{semi-cross} $\mathcal{B}(n,1,k,0)$ were investigated in \cite{stein_factoring_1967,szabo_problems_1986,hickerson_abelian_1986,szabo_factoring_2009}.  Schwartz \cite{schwartz_quasi-cross_2012} introduced the \emph{quasi-cross} $\mathcal{B}(n,1,k_1,k_2)$ with applications to flash memory.
	
	On the coding-theoretic side of the single-error case, systematic constructions of codes correcting asymmetric limited magnitude errors were given by Kl{\o}ve, Luo, Naydenova and Yari~\cite{klove2011asymmetric}, and perfect codes for asymmetric limited magnitude errors were studied by Buzaglo and Etzion~\cite{buzaglo2013chairs} via tilings of~$\Z^n$ with $n$-dimensional chairs; nonexistence results for lattice tilings by	quasi-crosses were subsequently obtained by	Schwartz~\cite{schwartz2014quasicross}.
	
	For $t\ge 2$, Zhang, Lian and Ge \cite{zhang_lattice_2023} recently gave a complete classification of lattice tilings by $\mathcal{B}(n,2,1,1)$.  Zhang and Ge \cite{zhang_new_2016} also studied codes correcting a single limited magnitude error and provided constructions for perfect and quasi-perfect codes.  In a different direction, the group ring method was employed by Leung, Tao, Wang and Zhang \cite{leung_lattice_2026} to make decisive progress on the $t=2$ asymmetric case $k_1>k_2$, establishing:
	\begin{enumerate}[label=(\roman*)]
		\item For $n\ge3$, a lattice tiling by $\mathcal{B}(n,2,3,0)$ exists iff $n=3$;
		\item No lattice tiling by $\mathcal{B}(n,2,k,k-1)$ exists for $n\ge3$, $k\ge2$;
		\item When $k_1+k_2+1$ is composite, no lattice tiling by $\mathcal{B}(n,2,k_1,k_2)$ exists for sufficiently large $n$.
	\end{enumerate}

	For general asymmetric balls $\mathcal{B}(n,t,k_1,k_2)$ with $k_1\neq k_2$, a systematic study of tilings was initiated by Wei and Schwartz~\cite{wei2022tilings}, who gave constructions based on perfect codes in the Hamming metric, proved several nonexistence results, and obtained complete classifications of lattice tilings in two special cases; 
	lattice packings and coverings by asymmetric limited magnitude balls were treated by Wei, Wang and Schwartz~\cite{wei2021packings}, and perfect codes correcting a single burst of limited magnitude errors were constructed by Wei and Schwartz~\cite{wei2022burst}. 
	
	The nonexistence results discussed above, together with those obtained in this paper, are summarized in Table~\ref{tab:alreadyexcluded}: the rows above the middle rule are due to previous work, while the rows below it record the results proved in the present paper.
	
	\begin{table}[t]
		\footnotesize
		\centering
		\caption{Summary of known nonexistence results for lattice tilings of
			$\Z^n$ by limited magnitude error balls $\mathcal{B}(n,t,k_1,k_2)$.}
		\label{tab:alreadyexcluded}
		\small
		\begin{tabular}{ccll}
			\toprule
			$t$ & $(k_1,k_2)$ & condition on $n$ for nonexistence & reference \\
			\midrule
			$1$ & $(k_1,k_2)$, quasi-cross
			& partial nonexistence results for lattice tilings
			& \cite{schwartz2014quasicross} \\
			$1$ & $(k,k)$, $(k,0)$
			& \makecell[l]{complete classifications of lattice tilings\\\quad by crosses and semi-crosses}
			& \cite{stein_factoring_1967,hickerson_abelian_1986,szabo_problems_1986} \\
			$2$ & $(1,1)$
			& complete classification of lattice tilings
			& \cite{zhang_lattice_2023} \\
			$2$ & $(3,0)$
			& for $n\ge3$: $n\neq3$ (tiling exists iff $n=3$)
			& \cite{leung_lattice_2026} \\
			$2$ & $(k,k-1)$, $k\ge2$
			& $n\ge3$
			& \cite{leung_lattice_2026} \\
			$2$ & $(k_1,k_2)$, $k_1>k_2$, $k_1+k_2+1$ composite
			& for sufficiently large $n$
			& \cite{leung_lattice_2026} \\
			\midrule
			$\ge1$ & $(k,k)$, $k\equiv4\pmod9$
			& \makecell[l]{if the base-$3$ addition of $t$ and $n-1-t$\\\quad produces exactly one carry}
			& Theorem~\ref{thm:keq4mod9} \\
			$n-1$ & $(k,k)$, $k\equiv8\pmod9$
			& if $n\equiv2,4\pmod6$
			& Theorem~\ref{thm:keq8mod9} \\
			$2$--$10$ & $ (k,k) $, $1\le k\le10$
			& \makecell[l]{residue classes of $n$ modulo $p^{2+\lfloor\log_p t\rfloor}$ for\\\quad primes $3\le p\le29$; shown in Appendix~\ref{app:computing-results}}
			& Theorem~\ref{thm:ideal} \\
			\bottomrule
		\end{tabular}
	\end{table}
	
	A complementary algebraic approach was formalized by Qureshi \cite{qureshi_non-existence_2020}, who introduced the \emph{Q-polynomial} and used it to convert the lattice tiling condition into a system of polynomial equations over finite fields. The definition of Q-polynomial is stated in Definition~\ref{def:Qpoly}.  This method was originally developed for Lee balls and yielded nonexistence results on linear perfect Lee codes in infinitely many dimensions.
	
	In detail, one can show that when $B$ is centrally symmetric ($B=-B$), $Q_B^{r}(X)$ is a homogeneous symmetric polynomial, hence expressible in the power sums $s_{2j}(X)$. That is,
	\begin{align}\label{eq:defclambda}
		Q_B^{r}(X)=\sum_{\lambda\in\mathcal{P}(r)}c_\lambda\prod_{i=1}^{\ell(\lambda)}s_{2\lambda_i}(X).
	\end{align}
	Moreover, the $c_\lambda$ are integers.
	
	Specifically, in the special case of Lee balls, Qureshi obtained an explicit formula only for \(c_{(r)}\), that is, the coefficient of the power sum \(s_{2r}\) in the expansion \eqref{eq:defclambda}. From these coefficients he derived a sufficient condition for the nonexistence of a lattice tiling, known as the \(p\)-condition, which depends solely on the \(c_{(r)}\)'s. 
	
	Furthermore, in the present author's joint preprint with Zhou \cite{zhou2026preprint},
	the complete expansion of the \(Q\)-polynomial for Lee balls was obtained, giving explicit formulas for all coefficients \(c_\lambda\), including the mixed ones with \(\ell(\lambda)>1\). Using this additional information, a new sufficient condition for the nonexistence of a lattice tiling was derived, which is formulated in terms of polynomial ideals. Moreover, it can be shown that this ideal-theoretic condition is strictly stronger than the \(p\)-condition.

	\subsection*{Our contribution}
	In this paper we focus on the symmetric case $k_1=k_2=k$.  For $X=(X_1,\dots,X_n)$,  denote the power sum symmetric polynomial by $s_j(X):=\sum_{i=1}^{n}X_i^{j}$. For a partition $\mu=(\mu_1,\dots,\mu_\ell)\vdash r$, let $\ell(\mu)$ denote its length and $m_j(\mu)$ the multiplicity of part $j$.
	
	Using enumerative combinatorial tools such as generating functions and Bell polynomials, we compute Q-polynomials associated with $\mathcal{B}(n,t,k,k)$ for any $n,t,k\in\Z^+$.
	\begin{thm}\label{thm:main}
		For $k\ge1$, $n\ge t\ge1$, and $r\ge1$, the $r$-th Q-polynomial of the limited magnitude error ball $\mathcal{B}(n,t,k,k)$ is
		\begin{align}\label{eq:Q-main}
			Q_{\mathcal{B}(n,t,k,k)}^{r}(X)=(2r)!\cdot[q^t]\Bigl(\frac{(1+2qk)^{n}}{1-q}\sum_{\mu\in\mathcal{P}(r)}\frac{1}{\prod_{j=1}^{r}m_j(\mu)!}\prod_{i=1}^{\ell(\mu)}d_{2\mu_i}(q,k)\,s_{2\mu_i}(X)\Bigr),
		\end{align}
		where $[q^t]f(q)$ denotes the coefficient of $q^t$ in the formal power series expansion of $f(q)$ and $d_{2j}(q,k)$ is defined by the following series:
		\[
		f(x):=\ln\!\Bigl(1+2q\sum_{b=1}^{k}\cosh(bx)\Bigr)=\sum_{j=0}^{\infty}d_{2j}(q,k)\,x^{2j}.
		\]
	\end{thm}
	Here, we present the explicit expressions for $d_0(q,k)$, $ d_2(q,k) $, $ d_4(q,k) $ and $ d_6(q,k) $.
	\begin{align*}
		d_0(q,k) &= \ln(2qk+1),\\[2mm]
		d_2(q,k) &= \frac{q k (2 k^{2}+3 k+1)}{12 q k+6},\\[2mm]
		d_4(q,k) &= -\frac{(4 k^{3} q+9 k^{2} q-3 k^{2}+7 q k-3 k+1)(k+1)(2 k+1) q k}{360(2 q k+1)^{2}},\\[2mm]
		d_6(q,k) &= \frac{(k+1)(2 k+1)
			\left(
			\begin{aligned}
				&64 q^{2} k^{6}+282 q^{2} k^{5}-90 k^{5} q+490 q^{2} k^{4}-243 k^{4} q\\
				&\quad+384 q^{2} k^{3}+9 k^{4}-210 k^{3} q+124 q^{2} k^{2}+18 k^{3}\\
				&\quad-36 k^{2} q+33 q k-9 k+3
			\end{aligned}
			\right) q k}{45360(2 q k+1)^{3}}.
	\end{align*}

	Moreover, combining properties of the Q-polynomial, a group-theoretic necessary condition for lattice tilings (Lemma~\ref{lem:tiling-group}), and power-sum properties in finite fields (Lemma~\ref{lem:powersumFF}) yields the following theorem, formulated in terms of polynomial ideals.

	\begin{thm}\label{thm:ideal}
		Let $p>2$ be a prime, and assume $\mathcal{B}(n,t,k,k)\subseteq\Z^n$ be a limited magnitude error ball such that $|\mathcal{B}(n,t,k,k)|=pM$ with $p\nmid M$. In $\F_p[Y_1,\dots,Y_{\frac{p-1}{2}}]$, define, for $1\le r\le\frac{p-1}{2}$,
		\[
		f_r:=\sum_{\lambda\in\mathcal{P}(r)}c_\lambda\prod_{i=1}^{\ell(\lambda)}Y_{\lambda_i},
		\]
		where $c_\lambda$ is as in \eqref{eq:defclambda}. If
		\[
		f_{\frac{p-1}{2}}\in\bigl(f_1,f_2,\dots,f_{\frac{p-1}{2}-1}\bigr),
		\]
		then no lattice tiling of $\Z^n$ by $\mathcal{B}(n,t,k,k)$ exists, where $(f_1,\dots,f_s)$ denotes the ideal generated by the polynomials $f_1,\dots,f_s$.
	\end{thm}
	
	Applying Theorem~\ref{thm:ideal}, for given $t,k$ and prime $p$, one can obtain congruence conditions on the dimension $n$ that rule out the existence of lattice tilings.  We record three representative results here; more details are shown in Example~\ref{ex:t6k1p5}, Example~\ref{ex:t3k2p5} and Example~\ref{ex:t4k3p11}.
	
	\begin{cor}\label{cor:t6k1}
		There is no lattice tiling of $\Z^n$ by $\mathcal{B}(n,6,1,1)$ if
		\[
		n\equiv 38,\;63,\;88,\;113\pmod{125}.
		\]
	\end{cor}
	
	\begin{cor}\label{cor:t3k2}
		There is no lattice tiling of $\Z^n$ by $\mathcal{B}(n,3,2,2)$ if
		\[
		n\equiv 7,\;8,\;12,\;13,\;17,\;18,\;22,\;23\pmod{25}.
		\]
	\end{cor}
	\begin{cor}\label{cor:t4k3}
		There is no lattice tiling of $\Z^n$ by $\mathcal{B}(n,4,3,3)$ if
		\[
		n\equiv 8,\;19,\;30,\;41,\;52,\;63,\;74,\;85,\;96,\;107,\;118\pmod{121}.
		\]
	\end{cor}
	
	Table~\ref{tab:excluded} lists, for each pair $(t,k)$ with $2\le t\le10$ and $1\le k\le10$, the percentage of dimensions $n$ for which Theorem~\ref{thm:ideal} rules out a lattice tiling of $\Z^n$ by $\mathcal{B}(n,t,k,k)$. Each entry is the fraction of all integers that this method excludes outright. The weakest entry rules out $4.99\%$ of all dimensions, the strongest $(t,k)=(10,5)$ rules out $76.69\%$ and the average over the table is $35.65\%$. How each entry is computed is explained in Example~\ref{ex:table-data}, and the complete raw data are collected in Appendix~\ref{app:computing-results}.
	
	\begin{table}[htbp]
		\footnotesize
		\centering
		\caption{Fraction (in percent) of integers $n$ ruled out by Theorem~\ref{thm:ideal} over all primes $3\le p\le29$}\label{tab:excluded}
		\begin{tabular}{c|cccccccccc}
			\toprule
			$t$ $\backslash$ $k$ & 1 & 2 & 3 & 4 & 5 & 6 & 7 & 8 & 9 & 10 \\
			\midrule
			2 & 22.22\% & 16.00\% & 12.24\% & 44.44\% &  8.26\% &  7.10\% & 34.67\% & 13.34\% &  4.99\% & 31.75\% \\
			3 & 14.81\% & 32.00\% & 24.49\% & 22.22\% & 16.53\% & 27.93\% & 42.07\% & 30.83\% &  9.97\% & 35.68\% \\
			4 & 14.81\% & 51.85\% & 42.49\% & 25.93\% & 30.36\% & 25.44\% & 59.36\% & 28.96\% & 14.96\% & 46.11\% \\
			5 & 18.52\% & 25.72\% & 55.51\% & 29.63\% & 42.98\% & 30.69\% & 28.95\% & 39.45\% & 24.38\% & 58.43\% \\
			6 & 41.08\% & 12.80\% & 61.22\% & 46.22\% & 41.32\% & 35.50\% & 45.10\% & 30.00\% & 27.33\% & 75.59\% \\
			7 & 34.13\% & 25.94\% & 13.36\% & 36.44\% & 49.59\% & 46.28\% & 48.75\% & 42.12\% & 44.54\% & 39.11\% \\
			8 & 20.27\% & 31.62\% & 16.22\% & 22.99\% & 58.59\% & 53.75\% & 40.71\% & 42.67\% & 41.15\% & 23.76\% \\
			9 & 28.73\% & 41.33\% & 26.73\% & 32.18\% & 66.12\% & 66.66\% & 57.10\% & 61.52\% & 48.50\% & 34.96\% \\
			10 & 22.64\% & 35.23\% & 32.05\% & 33.32\% & 76.69\% & 65.80\% & 39.35\% & 58.01\% & 46.80\% & 40.73\% \\
			\bottomrule
		\end{tabular}
	\end{table}
	
	Finally, we obtain two infinite families of parameters for which no lattice tiling exists. We specialize Corollary~\ref{cor:p=3}, which can be proved by fixing $p=3$ in Theorem~\ref{thm:conditionY}, to the congruence $3\mid[q^t]\frac{(1+2qk)^n}{1-q}B_2(0,2d_2(q,k))$; since $d_2(q,k)=\frac{qk(k+1)(2k+1)}{6(1+2qk)}$, this congruence is automatic for $k\equiv0,4,8\pmod9$. A case-by-case analysis of $|\mathcal{B}(n,t,k,k)|\bmod9$ then shows: for $k\equiv4\pmod9$ the obstruction becomes $v_3\binom{n-1}{t}=1$, which by Kummer's theorem means that the base-$3$ addition of $t$ and $n-1-t$ produces exactly one carry, while for $k\equiv8\pmod9$ and $t=n-1$ it reduces to $n\equiv2,4\pmod6$. This yields the following two theorems.
	\begin{thm}\label{thm:keq4mod9}
		There is no lattice tiling of $\Z^n$ by $\mathcal{B}(n,t,k,k)$ if $k\equiv4\pmod9$ and the addition of $t$ and $n-1-t$ in base $3$ produces exactly one carry.
	\end{thm}
	\begin{thm}\label{thm:keq8mod9}
		There is no lattice tiling of $\Z^n$ by $\mathcal{B}(n,n-1,k,k)$ if $k\equiv8\pmod9$ and $n\equiv2,4\pmod{6}$.
	\end{thm}
	
	\subsection*{Organization of the paper}
	The rest of the paper is organized as follows. Section~\ref{sec:pre} collects the preliminaries: power sums, integer partitions and complete Bell polynomials; the finite-field facts; the group characterization of lattice tilings; and Kummer's theorem. 
	In Section~\ref{sec:Qpoly} we prove Theorem~\ref{thm:main}, the explicit formula for the Q-polynomials of $\mathcal{B}(n,t,k,k)$.
	In Section~\ref{sec:ideal}, Theorem~\ref{thm:ideal} is proved through the chain Lemma~\ref{lem:conditionX} $\Rightarrow$ Theorem~\ref{thm:conditionY} $\Rightarrow$ Theorem~\ref{thm:ideal}, and illustrates the criterion with the examples $(t,k,p)=(6,1,5),(3,2,5),(4,3,11)$. Moreover, Example~\ref{ex:table-data} there explains how the entries of Table~\ref{tab:excluded} are computed. 
	In Section~\ref{sec:family} we obtain the two infinite families implying the nonexistence of lattice tilings by proving Theorems~\ref{thm:keq4mod9} and \ref{thm:keq8mod9}. Their sizes are estimated by numerical computations and an asymptotic analysis.
	The appendices contain the general pseudo-code for Theorems~\ref{thm:main} and \ref{thm:ideal}, and the complete raw data underlying Table~\ref{tab:excluded}. 
	
	\section{Preliminaries}\label{sec:pre}
	
	\subsection{Counting and combinatorial preliminaries}
	For $X=(X_1,\dots,X_n)$, let $s_j(X):=\sum_{i=1}^{n}X_i^{\,j}$ denote the power sum symmetric polynomial (with $s_0(X)=n$).
	
	Let $\mathcal{P}(r)$ be the set of integer partitions of $r$.  For $\lambda=(\lambda_1,\dots,\lambda_{\ell})\vdash r$, denote by $\ell(\lambda)$ the length and by $m_j(\lambda)$ the multiplicity of part $j$.
	
	The complete Bell polynomials $B_{r}(a_1,\dots,a_{r})$ are defined by
	\[
	\exp\!\Bigl(\sum_{j=1}^{\infty}a_j\frac{z^j}{j!}\Bigr)=\sum_{r=0}^{\infty}B_{r}(a_1,\dots,a_{r})\frac{z^{r}}{r!},
	\]
	with $B_0=1$.  See \cite{stanley_enumerative_2011} for the explicit formula
	\begin{align}\label{eq:Bell-explicit}
		B_{r}(a_1,\dots,a_{r})=\sum_{\lambda\in\mathcal{P}(r)}\frac{r!}{\prod_{i=1}^{\ell(\lambda)}\lambda_i!}\cdot\frac{1}{\prod_{j=1}^{r}m_j(\lambda)!}\cdot\prod_{i=1}^{\ell(\lambda)}a_{\lambda_i}.
	\end{align}
	
	Next, the Q-polynomial is defined as follows.
	\begin{defn}\label{def:Qpoly}
		For a finite set $B\subset\Z^n$ and $r\in\Z^+$, the \emph{Q-polynomial} of order $r$ is
		\[
		Q_B^{r}(X):=\sum_{\mathbf{b}\in B}\langle X,\mathbf{b}\rangle^{2r},\qquad X=(X_1,\dots,X_n),
		\]
		where $\langle,\rangle$ denotes the standard inner product.
	\end{defn}
	
	The size of the symmetric error ball \(\mathcal{B}(n,t,k,k)\) is
	\begin{align*}
		|\mathcal{B}(n,t,k,k)|=[q^t]\frac{(1+2qk)^n}{1-q}=\sum_{i=0}^{t}(2k)^i\binom{n}{i},
	\end{align*}
	where $[q^t]f(q)$ denotes the coefficient of $q^t$ in the formal power series expansion of $f(q)$.
	
	\subsection{Group-theoretic and finite-field preliminaries}
	The following lemma bridges lattice tilings with finite abelian groups \cite{horak_diameter_2012}.
	
	\begin{lem}\label{lem:tiling-group}
		Let $B$ be a finite subset of $\Z^n$ with $|B|=N$.  There exists a lattice tiling of $\Z^n$ by translates of $B$ if and only if there exists a finite abelian group $G$ of order $N$ and a group homomorphism $\phi:\Z^n\to G$ such that the restriction $\phi|_B:B\to G$ is a bijection.
	\end{lem}
	
	By the fundamental theorem of finite abelian groups, the following lemma holds.
	\begin{lem}\label{lem:mto1}
		Let $G$ be a finite abelian group of order $pM$, where $p$ is a prime. Then there exists a surjective homomorphism
		$\Phi:G\rightarrow\F_p$. In particular, $\Phi$ is an $M$-to-$1$ map.
	\end{lem}
	
	It is well known that power sums of elements in a finite field can be written as follows.
	\begin{lem}\label{lem:powersumFF}
		For a prime $p$, denote the finite field of characteristic $p$ by $\F_p$. For every integer $k\ge1$, we have
		\[
		\sum_{\alpha\in\F_p}\alpha^k=
		\begin{cases}
			-1,&\text{if }p-1\mid k;\\
			0,&\text{if }p-1\nmid k.
		\end{cases}
		\]
	\end{lem}
	
	\subsection{\(p\)-adic valuation and periodicity preliminaries}
	
	For the periodicity statement, we use the following uniform continuity result of Cahen--Chabert \cite{cahen1997integer}. 
	Let \(V\) be a rank-one	discrete valuation domain with finite residue field, and let \(v\) denote its valuation. 
	Let \(\mathfrak m\) be the maximal ideal of \(V\), and put \(q=\lvert V/\mathfrak m\rvert\); since the residue field is finite, \(q\) is a prime power. 
	Let \(\operatorname{Int}(V)\) denote the ring of integer-valued polynomials on \(V\).
	
	\begin{lem}[Cahen--Chabert \cite{cahen1997integer}, \S III.3, Proposition~III.3.1; see also \cite{gilmer1983finitely} for the case $V=\Z$]
		\label{lem:cahen-chabert}
		Let $h$ be an integer, and let $f\in\operatorname{Int}(V)$ satisfy
		$\deg(f)<q^{h}$. Then, for every integer $s$ and all $a,b\in V$,
		\[
		v(b-a)\ge s+h
		\quad\Longrightarrow\quad
		v\bigl(f(b)-f(a)\bigr)\ge s+1 .
		\]
	\end{lem}
	
	By substituting $V=\Z_p$ (so that $v=v_p$ and $q=p$),
	$h=\lfloor\log_p d\rfloor+1$ (so that $d<p^h$), $s=e-1$, $a=n$ and
	$b=n+p^{\,e+h-1}$ (note that $v_p(b-a)=e+h-1\ge(e-1)+h$), and viewing
	$f\in\operatorname{Int}(\Z)\subseteq\operatorname{Int}(\Z_p)$,
	Lemma~\ref{lem:cahen-chabert} yields $v_p\bigl(f(n+p^{\,e+h-1})-f(n)\bigr)\ge e$,
	that is, the following.
	
	\begin{lem}\label{lem:period}
		Let $f$ be an integer-valued polynomial of degree $d\ge 1$ and let $p$ be a prime.
		Then for every integer $n$ and every $e\ge 1$,
		\[
		f\bigl(n+p^{\,e+\lfloor\log_p d\rfloor}\bigr)\equiv f(n)\pmod{p^e}.
		\]
		In particular, the sequence $f(n)\bmod p^e$ is periodic in $n$ with period
		dividing $p^{\,e+\lfloor\log_p d\rfloor}$; if $d<p$, the period divides $p^e$.
	\end{lem}
	
	In the last section, the nonexistence criterion for the family $k\equiv4\pmod9$ (Theorem~\ref{thm:keq4mod9}) will be reduced to the condition $v_3\binom{n-1}{t}=1$. The following classical theorem of Kummer expresses the $p$-adic valuation of a binomial coefficient in terms of carries in base-$p$ addition, which makes such conditions easy to verify.
	
	\begin{thm}[Kummer's theorem, \cite{kummer_ueber_1852}]\label{thm:kummer}
		Let $p$ be a prime and let $a,b$ be nonnegative integers. Then the $p$-adic valuation $v_p\binom{a+b}{a}$ equals the number of carries when $a$ and $b$ are added in base $p$. In particular, for $0\le u\le n$, the valuation $v_p\binom{n}{u}$ equals the number of carries when $u$ and $n-u$ are added in base $p$.
	\end{thm}

	\section{Q-polynomial associated with symmetric error balls}\label{sec:Qpoly}
	
	In this section we prove Theorem~\ref{thm:main}.  Fix $k\ge 1$.  Introduce the generating function
	\begin{align*}
		F_n(X,z,q):=\sum_{t\ge0}\sum_{\mathbf{b}\in\mathcal{B}(n,t,k,k)}q^te^{z\langle X,\mathbf{b}\rangle}
	\end{align*}
	By the definition of Q-polynomial, expanding by the multinomial theorem, the coefficient of $z^{2r}q^t$ in $F_n$ is $\frac{Q_{\mathcal{B}(n,t,k,k)}^{r}(X)}{(2r)!}$, which gives
	\begin{align}\label{eq:dFdz1}
		\frac{\partial^{2r}F_n}{\partial z^{2r}}\Big|_{z=0}=\sum_{t\ge0}Q_{\mathcal{B}(n,t,k,k)}^{r}(X)q^t.
	\end{align}
	We only need to characterize the following expression:
	\[[q^t]\frac{\partial^{2r}F_n}{\partial z^{2r}}\Big|_{z=0}.\]

	Define the weight-exact shell by
	\[
	\mathcal{S}(n,t,k,k):=\{\mathbf{a}=(a_1,\dots,a_n)\in\Z^n:a_i\in[-k,k],\ \text{wt}(\mathbf{a})= t\}.
	\]
	Then $\mathcal{B}(n,t,k,k)=\bigcup_{s=0}^{t}\mathcal{S}(n,s,k,k)$ is a disjoint union, and
	\[
	Q_{\mathcal{B}(n,t,k,k)}^{r}(X)=\sum_{s=0}^{t}Q_{\mathcal{S}(n,s,k,k)}^{r}(X).
	\]

	\begin{proof}[Proof of Theorem~\ref{thm:main}]
		For a condition \(P\), let \(\mathbb{1}_{\{P\}}\) denote the indicator function that equals \(1\) if \(P\) holds and \(0\) otherwise. In particular, \(\mathbb{1}_{\{b\neq0\}}\) is \(1\) when \(b\neq0\) and \(0\) when \(b=0\). By geometric properties of limited magnitude error balls, we obtain
		\begin{align*}
			F_n(X,z,q)=&\sum_{t\ge0}\sum_{s=0}^{t}\sum_{\mathbf{b}\in\mathcal{S}(n,s,k,k)}q^te^{z\langle X,\mathbf{b}\rangle}=\frac{1}{1-q}\sum_{t\ge0}\sum_{\mathbf{b}\in\mathcal{S}(n,t,k,k)}q^te^{z\langle X,\mathbf{b}\rangle}\\
			=&\frac{1}{1-q}\sum_{\mathbf{b}\in([-k,k]\bigcap\Z)^n}q^{\text{wt}(\mathbf{b})}e^{z\langle X,\mathbf{b}\rangle}\\
			=&\frac{1}{1-q}\prod_{i=1}^{n}\sum_{b=-k}^{k}q^{\mathbb{1}_{\{b\neq0\}}}e^{zX_ib}\\
			=&\frac{1}{1-q}\prod_{i=1}^{n}\Bigl(1+q\sum_{b=1}^{k}(e^{zX_i b}+e^{-zX_i b})\Bigr)=\frac{1}{1-q}\prod_{i=1}^{n}\left(1+2q\sum_{b=1}^{k}\cosh(zX_i b)\right)
		\end{align*}
		
		Recall that $d_{2j}(q,k)$ is defined by the following series:
		\[
		f(x):=\ln\!\Bigl(1+2q\sum_{b=1}^{k}\cosh(bx)\Bigr)=\sum_{j=0}^{\infty}d_{2j}(q,k)\,x^{2j}.
		\]
		
		On the other hand, taking the logarithm gives
		\[
		\begin{aligned}
			F_n(X,z,q)
			&=\frac{1}{1-q}\prod_{i=1}^{n}\Bigl(1+q\sum_{b=1}^{k}(e^{zX_i b}+e^{-zX_i b})\Bigr)
			=\frac{1}{1-q}\exp\sum_{i=1}^{n}f(zX_i)  \\
			&=\frac{1}{1-q}\exp\sum_{j=0}^{\infty}d_{2j}(q,k)\,s_{2j}(X)\,z^{2j}
			=\frac{(1+2qk)^n}{1-q}\exp\!\Bigl(\sum_{j=1}^{\infty}d_{2j}(q,k)\,s_{2j}(X)\,(2j)!\,\frac{z^{2j}}{(2j)!}\Bigr).
		\end{aligned}
		\]
		Set $a_{2j-1}:=0$ and $a_{2j}:=d_{2j}(q,k)\,s_{2j}(X)\,(2j)!$ for $j\ge1$.  By the definition of complete Bell polynomials,
		\[
		F_n(X,z,q)=\frac{(1+2qk)^n}{1-q}\sum_{r=0}^{\infty}B_{r}(a_1,\dots,a_{r})\frac{z^{r}}{r!},
		\]
		and it follows that
		\[
		\frac{\partial^{2r}F_n}{\partial z^{2r}}\Big|_{z=0}
		=\frac{(1+2qk)^n}{1-q}\cdot B_{2r}(a_1,\dots,a_{2r}).
		\]
		Comparing \eqref{eq:dFdz1} with the expression below gives
		\[
		\sum_{t\ge0}Q_{\mathcal{B}(n,t,k,k)}^{r}(X)q^t=\frac{(1+2qk)^n}{1-q}\cdot B_{2r}(a_1,\dots,a_{2r}).
		\]
		Extracting the coefficient of $q^t$,
		\[
		Q_{\mathcal{B}(n,t,k,k)}^{r}(X)=[q^t]\Bigl(\frac{(1+2qk)^n}{1-q}\cdot B_{2r}(a_1,\dots,a_{2r})\Bigr).
		\]
		Finally, substitute the explicit formula~\eqref{eq:Bell-explicit} for $B_{2r}$.  Since only even parts contribute, partitions $\lambda\in\mathcal{P}(2r)$ correspond to $\mu\in\mathcal{P}(r)$ via $\lambda_i=2\mu_i$, giving
		\[
		B_{2r}(a_1,\dots,a_{2r})
		=\sum_{\mu\in\mathcal{P}(r)}\frac{(2r)!}{\prod_{j=1}^{r}m_j(\mu)!}\prod_{i=1}^{\ell(\mu)}d_{2\mu_i}(q,k)\,s_{2\mu_i}(X),
		\]
		which yields the claimed formula~\eqref{eq:Q-main}.
	\end{proof}

	\section{Ideal characterization of lattice tilings}\label{sec:ideal}
	In this section we prove Theorem~\ref{thm:ideal} and illustrate the resulting nonexistence criterion with three examples, from which we derive Corollary~\ref{cor:t6k1}, Corollary~\ref{cor:t3k2} and Corollary~\ref{cor:t4k3}.
	
	We prove Lemma \ref{lem:conditionX}, and then prove the following chain of implications:
	\begin{align}\label{eq:idealchain}
		\text{Lemma~\ref{lem:conditionX}}\Rightarrow\text{Theorem~\ref{thm:conditionY}}\Rightarrow\text{Theorem~\ref{thm:ideal}}.
	\end{align}
	\begin{lem}\label{lem:conditionX}
		Let $p>2$ be a prime. Assume that the limited magnitude error ball $\mathcal{B}(n,t,k,k)$ satisfies \(|\mathcal{B}(n,t,k,k)|=pM\) and there is no $X=(X_1,\dots,X_n)\in\F_p^n$ such that the following system of equations holds:
		\begin{align}
			\label{eq:systemX}
			\begin{cases}
				Q_{\mathcal{B}(n,t,k,k)}^{r}(X)\equiv0\pmod{p},\text{ for }1\le r<\frac{p-1}{2};\\
				Q_{\mathcal{B}(n,t,k,k)}^{\frac{p-1}{2}}(X)\equiv-M\pmod{p}.
			\end{cases}
		\end{align}
		Then there is no lattice tiling of \(\Z^n\) by $\mathcal{B}(n,t,k,k)$.
	\end{lem}
	\begin{proof}[Proof of Lemma \ref{lem:conditionX}]
		We prove the contrapositive. Suppose that there exists a lattice tiling of \(\Z^n\) by $\mathcal{B}(n,t,k,k)$. By Lemma \ref{lem:tiling-group}, there is an abelian additive group \(G\) and a homomorphism \(\varphi:\Z^n\rightarrow G\) such that the restriction \(\varphi|_{\mathcal{B}(n,t,k,k)}:\mathcal{B}(n,t,k,k)\rightarrow G\) is a bijection.
		Because $p\mid|\mathcal{B}(n,t,k,k)|$, Lemma~\ref{lem:mto1} gives a surjective homomorphism $\Phi:G\rightarrow(\F_p,+)$. Set $x_i=(\Phi\circ\varphi)(e_i)$, where $e_1,\dots,e_n$ are the standard basis vectors of $\Z^n$. Then
		\[
		(\Phi\circ\varphi)(\mathcal{B}(n,t,k,k))
		=\{*~a_1x_1+\cdots+a_nx_n:
		a_i\in\Z,\ a_i\in[-k,k],\ \text{wt}(\mathbf{a})\le t~*\}
		=M\F_p.
		\]
		where $\{*~\cdots~*\}$ denotes a multiset.
		
		By Lemma~\ref{lem:powersumFF}, we obtain
		\[
		\sum_{\alpha\in\Phi\circ\varphi(\mathcal{B}(n,t,k,k))}\alpha^{2r}=
		\begin{cases}
			0,&\text{if }p-1\nmid 2r;\\
			-M,&\text{if }p-1\mid 2r.
		\end{cases}
		\]
		
		Meanwhile, notice that \(\Phi\circ\varphi(\mathcal{B}(n,t,k,k))=\{*~\left\langle x,\mathbf{b}\right\rangle:\mathbf{b}\in \mathcal{B}(n,t,k,k)~*\}\), where \(x=(x_1,\dots,x_n)\in\F_p^n\), i.e.,
		\[
		\sum_{\mathbf{b}\in \mathcal{B}(n,t,k,k)}\left\langle x,\mathbf{b}\right\rangle^{2r}=
		\begin{cases}
			0,&\text{if }p-1\nmid 2r;\\
			-M,&\text{if }p-1\mid 2r.
		\end{cases}
		\]
		
		Since \(Q_{\mathcal{B}(n,t,k,k)}^{r}(X)=\sum_{\mathbf{b}\in \mathcal{B}(n,t,k,k)}\left\langle X,\mathbf{b}\right\rangle^{2r}\), \(x=(x_1,\dots,x_n)\) is a solution to the system \eqref{eq:systemX}.
	\end{proof}
	\begin{thm}\label{thm:conditionY}
		Let \(p>2\) be a prime. Assume that the limited magnitude error ball $\mathcal{B}(n,t,k,k)$ satisfies $|\mathcal{B}(n,t,k,k)|=pM$ and there do not exist \(Y_1,\dots,Y_{\frac{p-1}{2}}\in\F_p\) such that the following system of equations holds:
		\begin{align}
			\label{eq:systemY}
			\begin{cases}
				&\sum_{\lambda\in\mathcal{P}(r)}c_\lambda Y_{\lambda}\equiv0\pmod{p},\text{ for }1\le r<\frac{p-1}{2};\\
				&\sum_{\lambda\in\mathcal{P}(\frac{p-1}{2})}c_\lambda Y_{\lambda}\equiv-M\pmod{p},
			\end{cases}
		\end{align}
		where \(Y_{\lambda}=Y_{\lambda_1}\cdots Y_{\lambda_\ell}\) for \(\lambda=(\lambda_1,\dots,\lambda_\ell)\) and $c_\lambda$ is defined in \eqref{eq:defclambda}. Then there is no lattice tiling of \(\Z^n\) by $\mathcal{B}(n,t,k,k)$.
	\end{thm}
	\begin{proof}[Proof of the implication chain \eqref{eq:idealchain}]
		For the first implication, we prove the contrapositive. Suppose that a lattice tiling exists. Lemma~\ref{lem:conditionX} gives a solution $x=(x_1,\dots,x_n)\in\F_p^n$ of \eqref{eq:systemX}. Since
		\[
		Q_{\mathcal{B}(n,t,k,k)}^{r}(X)=\sum_{\lambda\in\mathcal P(r)}c_\lambda \prod_{i=1}^{\ell(\lambda)}s_{2\lambda_i}(X),
		\]
		the tuple
		\[
		(Y_1,\dots,Y_{(p-1)/2})
		=\left(\sum_{i=1}^n x_i^2,\sum_{i=1}^n x_i^4,\dots,
		\sum_{i=1}^n x_i^{p-1}\right)
		\]
		solves \eqref{eq:systemY}. This proves
		$\text{Lemma~\ref{lem:conditionX}}\Rightarrow
		\text{Theorem~\ref{thm:conditionY}}$.
		
		For the second implication, suppose that a lattice tiling exists and write
		$|\mathcal{B}(n,t,k,k)|=pM$ with $p\nmid M$. By the contrapositive of
		Theorem~\ref{thm:conditionY}, there is a solution
		$Y=(Y_1,\dots,Y_{(p-1)/2})$ of \eqref{eq:systemY}.
		Let
		\[
		I=\left(
		\sum_{\lambda\in\mathcal P(1)}c_\lambda \prod_{i=1}^{\ell(\lambda)}Y_{\lambda_i},
		\dots,
		\sum_{\lambda\in\mathcal P((p-3)/2)}c_\lambda \prod_{i=1}^{\ell(\lambda)}Y_{\lambda_i}
		\right).
		\]
		Every generator of $I$ vanishes at $Y$. Therefore, if
		$\sum_{\lambda\in\mathcal P((p-1)/2)}c_\lambda \prod_{i=1}^{\ell(\lambda)}Y_{\lambda_i}\in I$,
		that polynomial also vanishes at $Y$, contradicting the last equation in
		\eqref{eq:systemY}, because $-M\ne0$ in $\F_p$.
		Thus existence of a tiling implies the required nonmembership, which is
		the contrapositive of Theorem~\ref{thm:ideal}.
	\end{proof}

	Finally, we introduce three examples to show how Theorem~\ref{thm:ideal} works.
	\begin{example}[$k=1$, $t=6$, $p=5$]\label{ex:t6k1p5}
		Consider $\mathcal{B}(n,6,1,1)$, whose size is
		\[
		|\mathcal{B}(n,6,1,1)|=\sum_{i=0}^{6}2^{i}\binom{n}{i}
		=\frac{4}{45}n^{6}-\frac{16}{15}n^{5}+\frac{50}{9}n^{4}-\frac{40}{3}n^{3}
		+\frac{736}{45}n^{2}-\frac{28}{5}n+1.
		\]
		Let $p=5$.  The system~\eqref{eq:systemX} consists of two equations:
		\[
		Q_{\mathcal{B}}^{1}(X)=A(n)\,s_{2}(X)\equiv0\pmod{5},\qquad Q_{\mathcal{B}}^{2}(X)=B(n)\,s_{4}(X)+C(n)\,s_{2}^{2}(X)\equiv -M\pmod{5},
		\]
		where $M=|\mathcal{B}(n,6,1,1)|/5$.  Applying Theorem~\ref{thm:main} together with the Bell polynomial expansion gives
		\begin{align*}
			A(n)&=\frac{8}{15}n^{5}-\frac{20}{3}n^{4}+\frac{104}{3}n^{3}
			-\frac{256}{3}n^{2}+\frac{504}{5}n-42,\\[4pt]
			B(n)&=\frac{8}{15}n^{5}-\frac{44}{3}n^{4}+\frac{392}{3}n^{3}
			-\frac{1600}{3}n^{2}+\frac{5064}{5}n-726,\\[4pt]
			C(n)&=8n^{4}-96n^{3}+448n^{2}-912n+684.
		\end{align*}
		By Theorem~\ref{thm:ideal}, there do not exist lattice tilings by $\mathcal{B}(n,6,1,1)$ if $5\parallel |\mathcal{B}(n,6,1,1)|$ and $B(n)Y_2+C(n)Y_1^2\in(A(n)Y_1)\subseteq\F_{5}[Y_1,Y_2]$, where \(a\parallel b\) means that \(a\) divides \(b\) exactly, i.e. \(a\mid b\) but \(a^2\nmid b\).  
		
		Since $|\mathcal{B}(n,6,1,1)|=\sum_{i=0}^{6}2^i\binom{n}{i}$ is an integer-valued polynomial of degree $6$ in $n$, and the coefficients $A(n),B(n),C(n)$ are integer-valued polynomials of degree at most $5$, Lemma~\ref{lem:period} shows that the condition $5\parallel|\mathcal{B}(n,6,1,1)|$ (which involves $|\mathcal{B}(n,6,1,1)|\bmod 25$) has period dividing $5^{\,2+\lfloor\log_5 6\rfloor}=5^3=125$, while the ideal-membership condition, which involves only $A(n),B(n),C(n)\bmod 5$, has period dividing $5^{\,1+\lfloor\log_5 5\rfloor}=25$. 
		Hence all relevant congruences have period dividing $125$, and it suffices to verify the conditions for $n=0,1,\dots,124$; the remaining cases follow by periodicity. The periods $25$ and $121$ in Examples~\ref{ex:t3k2p5} and~\ref{ex:t4k3p11} follow the same way.

		This gives
		\[
		n\equiv 38,\;63,\;88,\;113\pmod{125},
		\]
		which establishes Corollary~\ref{cor:t6k1}.
	\end{example}
	
	\begin{example}[$k=2$, $t=3$, $p=5$]\label{ex:t3k2p5}
		Now consider $\mathcal{B}(n,3,2,2)$.  Its size is
		\[
		|\mathcal{B}(n,3,2,2)|=\sum_{i=0}^{3}4^{i}\binom{n}{i}
		=\frac{32}{3}n^{3}-24n^{2}+\frac{52}{3}n+1.
		\]
		Again take $p=5$.  The system~\eqref{eq:systemX} consists of two equations:
		\[
		Q_{\mathcal{B}}^{1}(X)=A(n)\,s_{2}(X)\equiv0\pmod{5},\qquad Q_{\mathcal{B}}^{2}(X)=B(n)\,s_{4}(X)+C(n)\,s_{2}^{2}(X)\equiv -M\pmod{5},
		\]
		where $M=|\mathcal{B}(n,3,2,2)|/5$.  Applying Theorem~\ref{thm:main} together with the Bell polynomial expansion gives
		\begin{align*}
			A(n)&=80n^{2}-200n+130,\\[4pt]
			B(n)&=272n^{2}-1880n+2542,\\[4pt]
			C(n)&=1200n-2100.
		\end{align*}
		Proceeding as in Example~\ref{ex:t6k1p5}, the Maple computation produces the excluded congruence classes
		\[
		n\equiv 7,\;8,\;12,\;13,\;17,\;18,\;22,\;23\pmod{25},
		\]
		which establishes Corollary~\ref{cor:t3k2}.
	\end{example}

	\begin{example}[$k=3$, $t=4$, $p=11$]\label{ex:t4k3p11}
		Now consider $\mathcal{B}(n,4,3,3)$.  Its size is
		\[
		|\mathcal{B}(n,4,3,3)|=\sum_{i=0}^{4}6^{i}\binom{n}{i}
		=54 n^{4}-288 n^{3}+504 n^{2}-264 n+1.
		\]
		Take $p=11$.  The system~\eqref{eq:systemX} consists of $ 5 $ equations:
		\begin{align*}
			Q_{\mathcal{B}}^{1}(X)=&A_1(n)\,s_{2}(X)\equiv0\pmod{11},\\
			Q_{\mathcal{B}}^{2}(X)=&A_2(n)\,s_{4}(X)+A_{1,1}(n)\,s_{2}^{2}(X)\equiv 0\pmod{11},\\
			\cdots\\
			Q_{\mathcal{B}}^{5}(X)=&\sum_{\lambda\in\mathcal{P}(5)}A_\lambda(n)\,\prod_{i=1}^{\ell(\lambda)}s_{2\lambda_i}(X)\equiv -M\pmod{11},
		\end{align*}
		where $M=|\mathcal{B}(n,4,3,3)|/11$.  Applying Theorem~\ref{thm:main} together with the Bell polynomial expansion, we obtain the following equations.
		
		\begin{align*}
			A_1(n)&=1008n^3 - 5544n^2 + 9744n - 5180,\\
			A_2(n)&=7056 n^{3}-81144 n^{2}+265776 n-264404,\\
			A_{1,1}(n)&=42336n^2 - 197568n + 228144,\\
			A_3(n)&=57168n^3 - 1796184n^2 + 11418864n - 19474340,\\
			A_{2,1}(n)&=1481760n^2 - 12841920n + 24778320,\\
			A_{1,1,1}(n)&=1975680n - 5597760,\\
			A_4(n)&=490896n^3 - 49311864n^2 + 609500976n - 1738103444,\\
			A_{3,1}(n)&=22409856n^2 - 491812608n + 1734236224,\\
			A_{2,2}(n)&=24202080n^2 - 306559680n + 872619440,\\
			A_{2,1,1}(n)&=193616640n - 935813760,\\
			A_{1,1,1,1}(n)&=64538880,\\
			A_5(n)&=4325328n^3 - 1509571224n^2 + 36717347184n - 173616608900,\\
			A_{4,1}(n)&=309264480n^2 - 21020247360n + 138453784560,\\
			A_{3,2}(n)&=1176517440n^2 - 30526231680n + 145040760480,\\
			A_{3,1,1}(n)&=4706069760n - 53993358720,\\
			A_{2,2,1}(n)&=10164873600n - 69459969600,\\
			A_{2,1,1,1}(n)&=13553164800,\\
			A_{1,1,1,1,1}(n)&=0.\\
		\end{align*}
		Proceeding as in Example~\ref{ex:t6k1p5}, the pseudo-code implementation in Appendix~\ref{app:maple-ideal} produces the excluded congruence classes
		\[
		n\equiv 8,\;19,\;30,\;41,\;52,\;63,\;74,\;85,\;96,\;107,\;118\pmod{121}
		\]
		which establishes Corollary~\ref{cor:t4k3}.
	\end{example}
	
	The origin of the data in Table~\ref{tab:excluded} is explained by the following example.
	
	\begin{example}[Derivation of the entry $(t,k)=(7,4)$ in Table~\ref{tab:excluded}]\label{ex:table-data}
		For each prime $3\le p\le29$, Theorem~\ref{thm:ideal} excludes a fixed set of residue classes of $n$ modulo $p^{2+\lfloor\log_p t\rfloor}$; here $t=7$, so the modulus is $p^3$ for $p=3,5,7$ and $p^2$ for $p\ge11$. The Maple computation (which uses Theorem~\ref{thm:main} to compute the Q-polynomials and then tests the ideal-membership condition for each residue class satisfying $p\mid|\mathcal{B}(n,t,k,k)|$ and $p^2\nmid|\mathcal{B}(n,t,k,k)|$) gives the following nonempty excluded sets:
		\begin{align*}
			p=3&:\ n\equiv 11,12,14,15,20,21,23,24\pmod{27};\\
			p=5&:\ n\equiv 16,17,41,42,67,91,116,117\pmod{125};\\
			p=7&:\ n\equiv 23,72,85,121,134,170,183,219,232,268,281,330\pmod{343};
		\end{align*}
		for $p=11,13,17,19,23,29$ the excluded set is empty. Thus the excluded fractions are $\frac{8}{27}$, $\frac{8}{125}$ and $\frac{12}{343}$. Since the moduli are pairwise coprime, the Chinese Remainder Theorem gives the total excluded fraction
		\[
		1-\Bigl(1-\frac{8}{27}\Bigr)\Bigl(1-\frac{8}{125}\Bigr)\Bigl(1-\frac{12}{343}\Bigr)
		=1-\frac{19\cdot117\cdot331}{27\cdot125\cdot343}
		\approx0.3644,
		\]
		which is exactly the entry $36.44\%$ of Table~\ref{tab:excluded}. Repeating this procedure for all $2\le t\le10$ and $1\le k\le10$ produces the whole table.
	\end{example}
	
	\section{Two Infinite Families of Parameters $(n,t,k)$ Implying No Lattice Tiling}\label{sec:family}
	In this section, by analyzing the coefficients of two $q$-series modulo $3$, we obtain two families of special parameters for which lattice tilings do not exist, shown as Theorem~\ref{thm:keq4mod9} and Theorem~\ref{thm:keq8mod9}.
	
	Recall that Theorem~\ref{thm:conditionY} gives a sufficient condition for the nonexistence of lattice tilings in terms of systems of multivariate equations. By fixing $p=3$ and assuming $p\nmid M$, it gives the following corollary.
	\begin{cor}\label{cor:p=3}
		There is no lattice tiling of $\Z^n$ by $\mathcal{B}(n,t,k,k)$ if 
		\[
		3\parallel|\mathcal{B}(n,t,k,k)| \text{ and }3\mid[q^t]\frac{(1+2qk)^n}{1-q}B_2(0,2d_2(q,k)),
		\]
		where \(B_{r}(a_1,\dots,a_{r})\) denotes the complete Bell polynomial; \([q^t]f(q)\) denotes the coefficient of \(q^t\) in the formal power series \(f(q)\); and
		\[
		d_2(q,k)=[x^2]\ln\!\left(1+2q\sum_{b=1}^{k}\cosh(bx)\right)
		=\frac{q\,k(k+1)(2k+1)}{6(1+2qk)}.
		\]
	\end{cor}
		Applying the preceding corollary, we compute
	\begin{align*}
		k\equiv0,4,8\pmod9&\Leftrightarrow3\mid\frac{k(k+1)(2k+1)}{3}\qquad(\text{direct check modulo }9)\\
		&\Rightarrow3\mid[q^t]\frac{q(1+2qk)^{n-1}}{1-q}\cdot\frac{k(k+1)(2k+1)}{3}\\
		&\Leftrightarrow3\mid[q^t]\frac{(1+2qk)^n}{1-q}B_2(0,2d_2(q,k)).
	\end{align*}
	
	Recall that 
	\[
	|\mathcal{B}(n,t,k,k)|=[q^t]\frac{(1+2qk)^n}{1-q}=\sum_{i=0}^{t}(2k)^i\binom{n}{i}.
	\]
	Next, we analyze the three residue classes \(k \equiv 0, 4, 8 \pmod 9\) separately and determine \(|\mathcal{B}(n,t,k,k)| \bmod 9\) in each case, which will yield Theorem~\ref{thm:keq4mod9} and Theorem~\ref{thm:keq8mod9}. 
	\begin{enumerate}
		\item $k\equiv0\pmod9$
		
		Then $|\mathcal{B}(n,t,k,k)|\equiv1\pmod9$, so the hypothesis $3\parallel|\mathcal{B}(n,t,k,k)|$ of Corollary~\ref{cor:p=3} is never satisfied.
		
		\item $k\equiv4\pmod9$
		
		\begin{proof}[Proof of Theorem~\ref{thm:keq4mod9}]
			In this case,
			\[
			|\mathcal{B}(n,t,k,k)|\equiv\sum_{i=0}^{t}(-1)^i\binom{n}{i}\pmod9\equiv(-1)^{t}\binom{n-1}{t}\pmod9,
			\]
			so that $3\parallel|\mathcal{B}(n,t,k,k)|\Leftrightarrow v_3\left(\binom{n-1}{t}\right)=1.$
			
			By Kummer's theorem (Theorem~\ref{thm:kummer}), the condition $v_3\left(\binom{n-1}{t}\right)=1$ is equivalent to the base-$3$ addition of $t$ and $n-1-t$ producing exactly one carry. Combined with Corollary~\ref{cor:p=3} --- whose second hypothesis $3\mid[q^t]\frac{(1+2qk)^n}{1-q}B_2(0,2d_2(q,k))$ holds by the computation above, since $k\equiv4\pmod9$ --- this proves Theorem~\ref{thm:keq4mod9}. For example, letting $n=3^a+1$ and $t=3^{a-1}$, we conclude that there is no lattice tiling of $\Z^n$ by $\mathcal{B}(n,t,k,k)$ for every integer $a\ge1$ and $k\equiv4\pmod9$.
		\end{proof}
		
		\item $k\equiv8\pmod9$
		
		\begin{proof}[Proof of Theorem~\ref{thm:keq8mod9}]
			In this case,
			\[
			|\mathcal{B}(n,t,k,k)|\equiv\sum_{i=0}^{t}(-2)^i\binom{n}{i}\pmod9.
			\]
			Unlike the previous case, the partial sum \(\sum_{i=0}^{t}(-2)^i\binom{n}{i}\) has no simple closed form for general \(t\). However, by imposing the restriction \(t=n-1\), we obtain
			\[
			|\mathcal{B}(n,n-1,k,k)|\equiv\sum_{i=0}^{n-1}(-2)^i\binom{n}{i}\pmod9\equiv(-1)^n-(-2)^n\pmod9.
			\]
			Since $(-2)^3\equiv1\pmod9$, it suffices to check $n$ modulo $6$; a direct computation gives
			\[
			v_3\bigl((-1)^n-(-2)^n\bigr)=1\Leftrightarrow n\equiv2,4\pmod6,
			\]
			which, combined with Corollary~\ref{cor:p=3} (whose second hypothesis holds since $k\equiv8\pmod9$), proves Theorem~\ref{thm:keq8mod9}.
		\end{proof}
	\end{enumerate}
	
	We close this section with a few words on the strength of Theorems~\ref{thm:keq4mod9} and \ref{thm:keq8mod9}. Theorem~\ref{thm:keq8mod9} rules out the ball $\mathcal{B}(n,n-1,k,k)$ for exactly two residue classes of $n$ modulo $6$, i.e. for one third of all dimensions, whenever $k\equiv8\pmod9$; since $k\equiv8\pmod9$ covers one ninth of all $k$, assuming independence, this family excludes $1/27\approx3.70\%$ of the pairs $(n,k)$ for the extremal weight $t=n-1$.
		
	We next investigate the strength of the following condition stated in Theorem~\ref{thm:keq4mod9}.
	\[
	v_3\!\left(\binom{n-1}{t}\right)=1.
	\]
	More precisely, we study the proportion of pairs \((n,t)\) in the set
	\[
	\{(n,t)\in[1,M]^2:\ n-1\ge t\}
	\]
	that satisfy this condition. Both numerical computations and an asymptotic analysis are presented.

	For \(M=100,500,1000\), a direct computation gives the following table. Here ``count'' denotes the number of pairs \((n,t)\in[1,M]^2\) with \(n-1\ge t\) and
	\[
	v_3\!\left(\binom{n-1}{t}\right)=1,
	\]
	while ``total'' is the total number of pairs in the triangular region,
	\[
	\#\{(n,t)\in[1,M]^2:\ n-1\ge t\}=\frac{M(M-1)}2.
	\]
	The ratio is the quotient of these two quantities.
	
	\begin{table}[h]
		\centering
		\begin{tabular}{cccc}
			\hline
			\(M\) & count & total & ratio \\
			\hline
			\(100\)  & \(1292\)  & \(4950\)   & \(0.261010\) \\
			\(500\)  & \(27546\) & \(124750\) & \(0.220810\) \\
			\(1000\) & \(90954\) & \(499500\) & \(0.182090\) \\
			\hline
		\end{tabular}
		\caption{Numerical values of the proportion of pairs satisfying \(v_3\!\left(\binom{n-1}{t}\right)=1\).}
	\end{table}
	
	As \(M\to\infty\), the strength of this condition is described by the following asymptotic theorem.

	\begin{thm}\label{thm:asymptotic-ratio}
		For every $M\ge 9$, let
		\[
		P(M):=
		\frac{
			\#\{(n,t)\in[1,M]^2:\ n-1\ge t,\ v_3\!\left(\binom{n-1}{t}\right)=1\}
		}{
			\#\{(n,t)\in[1,M]^2:\ n-1\ge t\}
		}.
		\]
		Then, 
		\[
		P(M)=\Theta(M^{\log_3 (2/3)}\ln M).
		\]
		In particular, $P(M)\to0$ as \(M\to\infty\).
	\end{thm}
	
	\begin{proof}
		The numerator equals
		\[
		E(M):=\sum_{u=0}^{M-1}\#\bigl\{0\le t\le u:\ v_3\!\tbinom{u}{t}=1\bigr\},
		\]
		since $v_3\tbinom{u}{0}=v_3(1)=0$. By Kummer's theorem, $v_3\tbinom{u}{t}$
		is the number of carries in the base-$3$ addition of $t$ and $u-t$; writing
		$s=u-t$,
		\[
		E(M)=\#\bigl\{(t,s)\in\Z_{\ge0}^{\,2}:\ t+s<M,\
		\text{$t+s$ has exactly one carry in base $3$}\bigr\},
		\]
		so that
		\[
		P(M)=\frac{E(M)}{\#\{(n,t)\in[1,M]^2:\ n-1\ge t\}}=\frac{E(M)}{M(M-1)/2}.
		\]
		
		For $k\ge1$ put
		\begin{align*}
			A_k&=\#\{(t,s)\in[0,3^k)^2:\ t+s<3^k,\ \text{exactly one carry}\},\\
			B_k&=\#\{(t,s)\in[0,3^k)^2:\ t+s\ge3^k,\ \text{exactly one carry}\}.
		\end{align*}
		For $ A_{k+1} $ and $ B_{k+1} $, write $t=\alpha3^k+t'$, $s=\beta3^k+s'$ with $\alpha,\beta\in\{0,1,2\}$,
		$0\le t',s'<3^k$, and let $c_k=\lfloor(t'+s')/3^k\rfloor\in\{0,1\}$ be the
		carry out of the low block; the total number of carries in $t+s$ is
		\[
		c(t,s)=c(t',s')+\Bigl\lfloor\frac{\alpha+\beta+c_k}{3}\Bigr\rfloor.
		\]
		If there is no top carry ($\alpha+\beta+c_k\le2$), then $c(t,s)=1$ iff
		$c(t',s')=1$: the pairs $(t',s')\in A_k$ (i.e.\ $c_k=0$) require
		$\alpha+\beta\le2$ ($6$ choices), the pairs $(t',s')\in B_k$ (i.e.\ $c_k=1$)
		require $\alpha+\beta\le1$ ($3$ choices). If the top carry is present
		($\alpha+\beta+c_k\ge3$), then $c(t,s)=c(t',s')+1$, so $c(t,s)=1$ iff the
		low block has no carry at all: there are
		$C_k=\#\{(t',s'):\ t'_i+s'_i\le2\ \forall i\}=6^k$ such pairs, and
		$\alpha+\beta\ge3$ allows $3$ choices. Hence
		\[
		A_{k+1}=6A_k+3B_k,\qquad B_{k+1}=3\cdot6^{\,k},
		\]
		with $A_1=0$, $B_1=3$; solving gives
		\[
		B_k=3\cdot6^{\,k-1},\qquad A_k=9(k-1)6^{\,k-2}\ (k\ge2),\qquad
		E(3^k)=A_k.
		\]
		
		Let $k=\lfloor\log_3 M\rfloor$, i.e., $3^k\le M<3^{k+1}$. By $M\ge 9$, we can know $k\ge2$.
		\begin{itemize}
			\item \textit{Upper bound.} Since $E$ is nondecreasing, $E(M)\le E(3^{k+1})
			=9k\,6^{k-1}$, while $M^{\log_3 (6)}\ge6^k$, $\ln M\ge k\ln 3$ and
			$\frac{M}{M-1}\le2$; therefore
			\[
			\frac{P(M)}{M^{\log_3 (2/3)}\ln M}
			=\frac{2E(M)}{M^{\log_3 (6)}\ln M}\,\frac{M}{M-1}
			\le\frac{2\cdot9k\,6^{k-1}}{6^k\,k\ln 3}\cdot2=\frac{6}{\ln 3}.
			\]
			\item \textit{Lower bound.} If $M\ge9$, then $k\ge2$; using $E(M)\ge E(3^k)$,
			$M^{\log_3 (6)}\le6^{k+1}$, $\ln M\le(k+1)\ln 3$ and
			$(k-1)/(k+1)\ge\frac13$,
			\[
			\frac{P(M)}{M^{\log_3 (2/3)}\ln M}
			\ge\frac{2E(M)}{M^{\log_3 (6)}\ln M}
			\ge\frac{18(k-1)6^{k-2}}{6^{k+1}(k+1)\ln 3}
			=\frac{k-1}{12(k+1)\ln 3}
			\ge\frac{1}{36\ln 3}.
			\]
		\end{itemize}

		Thus the two-sided bound holds for all $M\ge9$, giving
		$P(M)=\Theta(M^{\log_3 (2/3)}\ln M)$; since $\log_3 (2/3)<0$, $P(M)\to0$.
	\end{proof}
	\section{Concluding Remarks}\label{sec:conc}
	\begin{enumerate}[label=(\arabic*)]
		
		\item In Section~\ref{sec:family} we treated only the case $p=3$, where,
		after imposing a restriction on $k$, the divisibility
		\[
		3 \mid [q^t]\,\frac{(1+2qk)^n}{1-q}\,B_2\bigl(0,\,2d_2(q,k)\bigr)
		\]
		led to the two infinite families of Theorems~\ref{thm:keq4mod9} and~\ref{thm:keq8mod9}. It would be	interesting to study the expansions of the series involved more closely and to obtain analogous infinite families for $p=5,7$ and, ideally, for every odd	prime $p$; progress in this direction would provide evidence for the conjecture that for each pair $(t,k)$ there are infinitely many $n$ for which	no lattice tiling of $\Z^n$ by $\mathcal{B}(n,t,k,k)$ exists. A second natural idea is to use the Chinese remainder theorem to combine the excluded residue classes obtained for several primes---equivalently, to work with the $p$-primary decomposition of the group $G$ of Lemma~\ref{lem:tiling-group} rather than with a single surjection $\Phi\colon G\to\F_p$---thereby recovering part of
		the information lost when passing from $G$ to a finite field.
		
		\item The group ring method of Leung, Tao, Wang and Zhang \cite{leung_lattice_2026} made decisive progress on the asymmetric balls $\mathcal{B}(n,2,k_1,k_2)$ with $k_1>k_2$, notably proving that no lattice tiling by $\mathcal{B}(n,2,k,k-1)$ exists for $n\ge3$, $k\ge2$. Since $Q$-polynomials and group ring equations appear to encode the same information about tilings, it would be desirable to develop a framework---for instance a non-symmetric analogue of the generating function $F_n(X,z,q)$, based on the weights	$e^{z\langle X,b\rangle}$ with $b\in[-k_2,k_1]$ in place of the symmetric	$\cosh$ expansion---through which the $Q$-polynomial method can be extended	to the asymmetric case.
		
	\end{enumerate}
	
	\bigskip
	\paragraph*{\bf Declaration of competing interest}
	The authors declare that they have no known competing financial interests or personal 
	relationships that could have appeared to influence the work reported in this paper.
	\bigskip
	\paragraph*{\bf Declaration of generative AI use}
	The authors used a generative-AI tool solely to assist with English-language editing. The authors reviewed and approved all revisions and remain fully responsible for the content of the manuscript.
	\bigskip
	\paragraph*{\bf Formalized and machine-checked in Lean 4}
	All statements proved in this paper, except those that rest on our Maple computations, have been formalized and machine-checked in the proof	assistant Lean~4 \cite{moura_lean4_2021} with the library Mathlib	\cite{mathlib_2020}. The formalization was produced with the help of the AI system Aristotle \cite{aristotle_harmonic} by Harmonic. It contains no unproved steps and uses only the standard axioms of Lean. 
	
	The authors have checked that the formal statements agree with the statements in the paper. The source code is available at
	\href{https://github.com/Xiao-00000/Symmetric-Limited-Magnitude-Error-Ball/tree/main/lean4}{this link}.
	
	Specifically, Theorem~\ref{thm:main} (and its Bell-polynomial form), the closed forms of $d_0, d_2, d_4, d_6$, the integrality of $c_\lambda$ in \eqref{eq:defclambda}, Lemmas~\ref{lem:tiling-group}, \ref{lem:mto1}, \ref{lem:powersumFF}, \ref{lem:cahen-chabert}, \ref{lem:period}, \ref{lem:conditionX} (along with Kummer's theorem), Theorems~\ref{thm:conditionY}, \ref{thm:ideal}, \ref{thm:keq4mod9}, \ref{thm:keq8mod9}, \ref{thm:asymptotic-ratio}, Corollaries~\ref{cor:p=3}, \ref{cor:t6k1}, \ref{cor:t3k2}, \ref{cor:t4k3}, Examples~\ref{ex:table-data}, \ref{ex:t6k1p5}, \ref{ex:t3k2p5}, \ref{ex:t4k3p11}, the general test derived from Theorem~\ref{thm:ideal}, and the numerical table in Section~\ref{sec:family} are formally verified.

	\appendix
	
	\section{Maple Implementation of Theorem~\ref{thm:ideal}}
	\label{app:maple-ideal}
	
	The following pseudo-code implements Theorem~\ref{thm:ideal} for arbitrary $t$, $k$ and an odd prime $p$: it computes the Q-polynomials $Q_{\mathcal{B}(n,t,k,k)}^{r}(X)$ via Theorem~\ref{thm:main}, collects the residue classes of $n$ satisfying $p\mid|\mathcal{B}(n,t,k,k)|$ and $p^2\nmid|\mathcal{B}(n,t,k,k)|$, and tests the ideal-membership condition. The period can be obtained by Lemma~\ref{lem:period}, which equals $p^{2+\lfloor\log_p t\rfloor}$.
	
	\newcounter{alg}
	\begin{tcolorbox}[colback=white,colframe=black,boxrule=0.8pt]
		\refstepcounter{alg}
		\textbf{Algorithm \thealg:} \textsc{Lattice-Tiling Nonexistence Tester}
		\label{alg:nonexistence-tester}
		
		\textbf{Input: }weight bound $t\ge2$, magnitude $k\ge1$, and an odd prime $p$
		
		\textbf{Output: }the set of residue classes $n \bmod p^{\,2+\lfloor\log_p t\rfloor}$ for which Theorem~\ref{thm:ideal} rules out a lattice tiling of $\Z^n$ by $\mathcal{B}(n,t,k,k)$
		\vspace{2pt}
		
		$\mathrm{period}\gets p^{\,2+\lfloor\log_p t\rfloor}$ \hfill $\triangleright$ \textnormal{period of $|\mathcal{B}(n,t,k,k)|\bmod p^2$ (Lem.~\ref{lem:period})}
		
		$C\gets\emptyset$
		
		\textbf{for} $n=1,\dots,\mathrm{period}$ \textbf{do}
		
		\hspace*{2em}compute $|\mathcal{B}(n,t,k,k)|=\sum_{i=0}^{t}(2k)^i\binom{n}{i}$
		
		\hspace*{2em}\textbf{if} $p\mid|\mathcal{B}(n,t,k,k)|$ and $p^2\nmid|\mathcal{B}(n,t,k,k)|$ \textbf{then}
		
		\hspace*{2em}\hspace*{2em}$C\gets C\cup\{n\}$ \hfill $\triangleright$ \textnormal{$C$ contains $n$ satisfying hypothesis of Thm.~\ref{thm:ideal}}
		
		\textbf{end for}
		
		Determine $d_{2j}(q,k)$ from $\ln\!\Bigl(1+2q\sum_{b=1}^{k}\cosh(bx)\Bigr)=\sum_{j\ge0}d_{2j}(q,k)\,x^{2j}$
		
		\textbf{for} $r=1,\dots,(p-1)/2$ \textbf{do}
		
		\hspace*{2em}$a_{2j-1}\gets0$,\quad $a_{2j}\gets(2j)!\,d_{2j}(q,k)\,Y_j$ for $1\le j\le r$
		
		\hspace*{2em}$Q^r\gets[q^t]\,\dfrac{(1+2qk)^n}{1-q}\,B_{2r}(a_1,\dots,a_{2r})$ \hfill $\triangleright$ \textnormal{by Thm.~\ref{thm:main}, with $Y_j= s_{2j}(X)$}
		
		\textbf{end for}
		
		$E\gets\emptyset$
		
		\textbf{for} each $n$ in $C$ \textbf{do}
		
		\hspace*{2em}let $J=\bigl(Q^1(n),\dots,Q^{(p-1)/2-1}(n)\bigr)$ be an ideal of $\F_p[Y_1,\dots,Y_{(p-1)/2}]$
		
		\hspace*{2em}\textbf{if} $Q^{(p-1)/2}(n)\in J$ \textbf{then}
		
		\hspace*{2em}\hspace*{2em}$E\gets E\cup\{n\}$ \hfill $\triangleright$ \textnormal{no lattice tiling (Thm.~\ref{thm:ideal})}
		
		\textbf{end for}
		
		\textbf{return} $E$  
	\end{tcolorbox}

\section{Computing Results for $2\le t\le10$, $1\le k\le10$ and prime $3\le p\le29$}
\label{app:computing-results}

This appendix collects the raw data produced by the Maple implementation of Theorem~\ref{thm:ideal} given in Appendix~\ref{app:maple-ideal}. For every pair $(t,k)$ with $2\le t\le10$ and $1\le k\le10$, and every prime $3\le p\le29$, the program outputs the residue classes of $n$ modulo $p^{2+\lfloor\log_p t\rfloor}$ that satisfy $p\mid|\mathcal{B}(n,t,k,k)|$ and $p^2\nmid|\mathcal{B}(n,t,k,k)|$ and pass the ideal-membership test of Theorem~\ref{thm:ideal}; for each such class, Theorem~\ref{thm:ideal} guarantees that no lattice tiling of $\Z^n$ by $\mathcal{B}(n,t,k,k)$ exists. These residue classes are precisely the data from which the percentages in Table~\ref{tab:excluded} were computed. Since the complete output is too long to be included in print, it is hosted at \url{https://github.com/Xiao-00000/Symmetric-Limited-Magnitude-Error-Ball}.

\end{document}